\documentclass[a4paper,twoside,12pt]{article}

\usepackage{amsmath}
\usepackage{amsfonts}
\usepackage{amssymb}
\usepackage{amsthm}
\usepackage{amsopn}
\usepackage{amscd}

\usepackage{pstricks}
\usepackage[latin1]{inputenc}
\usepackage{courier}
\usepackage{url}

\usepackage[vcentering,dvips]{geometry}
\def\address#1#2{\begingroup
\noindent\parbox[t]{7.8cm}{%
\small{\scshape\ignorespaces#1}\par\vskip1ex
\noindent\small{\itshape E-mail address}%
\/: #2\par\vskip4ex}\hfill%
\endgroup}%

\author{Nathan Owen Ilten}
\title{Calculating Milnor Numbers and Versal Component Dimensions from P-Resolution Fans}
\date{}
\newcommand{\rsuchthat}{\ \right|\left.\ }
\newcommand{\lsuchthat}{\ \right.\left|\ }

\renewcommand{\{}{\left\lbrace\left.  }
\renewcommand{\}}{\right\rbrace\right.} 

\DeclareMathOperator{\spec}{Spec}

\DeclareMathOperator{\tv}{TV}

\DeclareMathOperator{\conv}{Conv}

\DeclareMathOperator{\qdef}{Def'}

\newtheorem{lemma}{Lemma}[section]
\newtheorem{prop}[lemma]{Proposition}
\newtheorem{cor}[lemma]{Corollary}
\newtheorem{thm}[lemma]{Theorem}

\theoremstyle{definition}
\newtheorem*{ex}{Example}
\newtheorem*{remark}{Remark}
\newtheorem*{defn}{Definition}

\begin{document}
\maketitle
\begin{abstract}
	We use Altmann's toric fan description of P-resolutions \cite{MR1652476} to formulate a new description of deformation theory invariants for two-dimensional cyclic quotient singularities. In particular, we show how to calculate the dimensions of the (reduced) versal base space components as well as Milnor numbers of smoothings over them. 
\end{abstract}

\section{Introduction}\label{intro}
The deformation theory of (two-dimensional) cyclic quotient  singularities is fairly well understood. Kollár and Shepherd-Barron have proven a correspondence between so-called P-resolutions and reduced components of the versal base space in \cite{MR922803}. P-Resolutions were further studied using continued fractions by Christophersen and Stevens in \cite{MR1129026} and \cite{MR1129040} respectively, who both managed to write down explicit equations for the reduced components of the versal deformation. In \cite{MR1652476}, Altmann uses the continued fractions of Christophersen and Stevens to describe P-resolutions in toric terms, that is, in terms of a fan.

This paper deals with the dimension of the versal base components as well as the Milnor numbers of smoothings over them. A formula for the Milnor numbers is provided in \cite{MR1129040}. Furthermore, for the $\mathbb{Q}$-Gorenstein one-parameter smoothing of $T$-singularities, Altmann has already provided a simple formula in toric terms. The first aim of this paper is to generalize Altmann's formula to smoothings over all components of any cyclic quotient singularity.  Our new formula allows the Milnor number of a smoothing to be read directly from the geometry of the fan describing the corresponding P-resolution.

A method for calculating the dimension of the versal base components was first provided by Kollár and Shepherd-Barron in \cite{MR922803}. Even better, the explicit equations in \cite{MR1129026} and \cite{MR1129040} allow one to write down a simple formula.  The second aim of this paper is to translate this component dimension formula into toric language using Altmann's toric description of P-resolutions. Our new  formula allows the dimension of a component to be read directly from the geometry of the fans describing the minimal resolution and the P-resolution corresponding to that component. Furthermore, the difference in dimension between two components can be read solely from the two fans describing the two corresponding components. Note also that our proofs of the two formulae, while not differing significantly from those by Stevens, can be easily understood within the context of toric geometry.  

In section \ref{prelim}, we provide necessary definitions and notation. We have chosen notation so as to be completely consistent with \cite{MR1652476}; in fact, readers familiar with this paper can probably skip section \ref{prelim}. In section \ref{milnor}, we describe Stevens' formula for Milnor numbers and  state and prove our new toric formula. Likewise, in section \ref{formulas} we present the existing component dimension formula and then state and prove our new toric formula. We finish in section \ref{example} by providing an example demonstrating the practicality of our formulae.

\section{Cyclic Quotients and P-Resolutions}\label{prelim}
In the following, we recall the notions of cyclic quotients and P-resolutions, as well as fixing notation. References are \cite{MR1234037} for toric varieties, and \cite{MR1652476} for P-resolutions.

Let $n$ and $q$ be relatively prime integers with $n\geq2$ and $0<q<n$. Let $\xi$ be a primitive $n$-th root of unity. The cyclic quotient singularity $Y_{(n,q)}$ is the quotient ${\mathbb{C}^2}/{(\mathbb{Z}/n\mathbb{Z})}$ where $\mathbb{Z}/n\mathbb{Z}$ acts on $\mathbb{C}^2$ via the matrix
	\begin{equation*}
		\left(\begin{array}{cc}\xi&0\\0&\xi^q\end{array}\right).
	\end{equation*}
	Every two-dimensional cyclic quotient singularity is in fact a two-dimensional toric variety: Let $N$ be a rank two lattice with dual lattice $M$ and let $\sigma\subset N\otimes \mathbb{R}$ be the cone generated by $(1,0)$ and $(-q,n)$. $Y_{(n,q)}$ is then isomorphic to the toric variety $U_\sigma=\spec \mathbb{C}[M\cap\sigma^\vee]$. Since by correct choice of basis every singular two-dimensional cone has generators $(1,0)$ and $(-q,n)$ for some $q$ and $n$ as above, every affine singular two-dimensional toric variety is a cyclic quotient singularity.

Introduced by Kollár and Shepherd-Barron in \cite{MR922803}, P-resolutions have proven key to understanding the deformation theory of cyclic quotient singularities.
\begin{defn}
Let $Y$ be a two-dimensional cyclic quotient singularity. A P-resolution of Y is a partial resolution $f:\widetilde{Y}\to Y$ containing only T-singularities such that the canonical divisor $K_{\widetilde{Y}}$ is ample relative to $f$. T-singularities are exactly those cyclic quotients admitting a $\mathbb{Q}$-Gorenstein one-parameter smoothing
\end{defn}
The following theorem describes the relationship between P-resolutions and reduced components of the base space of the versal deformation for a two-dimensional cyclic quotient singularity $Y$:
\begin{thm}
	There is a bijection between the set of P-resolutions $\{\widetilde{Y}_\nu\}$ of $Y$ and the components of the reduced versal base space $S$, induced by the natural maps $\qdef \widetilde{Y}_\nu\to S$, where $\qdef \widetilde{Y}_\nu$ is the space of $\mathbb{Q}$-Gorenstein deformations of $\widetilde{Y}_\nu$.
\begin{proof}
	See \cite{MR922803}, section 3.
\end{proof}
\end{thm}

P-Resolutions were described by Christophersen and Stevens in terms of continued fractions in \cite{MR1129026} and \cite{MR1129040}.  In \cite{MR1652476}, Altmann provides a toric description in terms of a fan; it is this latter description that we shall use in our dimension formula.

Let $c_1,c_2,\ldots,c_k\in\mathbb{Z}$. The continued fraction $[c_1,c_2,\ldots,c_k]$ is inductively defined as follows if no division by $0$ occurs: $[c_k]=c_k$,  $[c_1,c_2,\ldots,c_k]=c_1-1/[c_2,\ldots,c_k]$. Now, if one requires that $c_i\geq2$ for every coefficient, each continued fraction yields a unique rational number.

Let $n$ and $q$ be relatively prime integers with $n\geq3$ and $0<q<n-1$.\footnote{This restriction simply ensures that $Y_{(n,q)}$ isn't a hypersurface, in which case the versal base space is irreducible.}  We consider the cyclic quotient singularity $Y_{(n,q)}$. Let $[a_2,a_3,\ldots,a_{e-1}],\ a_i\geq2$ be the unique continued fraction expansion of $n/(n-q)$. Note that $e$ equals the embedding dimension of $Y_{(n,q)}$. Furthermore, the generators of the semigroup $M\cap\sigma^\vee{}$ are related to this continued fraction: $w^1=[0,1]$, $w^e=[n,q]$, and $w^{i-1}+w^{i+1}=a_i w^i$.

Likewise, let $[b_1,\ldots,b_r]$ be the unique continued fraction expansion of $n/q$.  The generators of the semigroup $N\cap\sigma$ are related to this continued fraction: $v^0=(1,0)$, $v^{r+1}=(-q,n)$, and $v^{i-1}+v^{i+1}=b_i v^i$. Drawing rays through the $v^i$ gives a polyhedral subdivision $\Sigma$ of $\sigma$. The corresponding toric variety $\tv(\Sigma)$ is the minimal resolution of $Y$ with self intersection numbers $-b_i$; the number of exceptional divisors in this resolution is $r$.

For a chain of integers $(k_2,\ldots,k_{e-1})$  define the sequence $q_1,\ldots,q_e$ inductively: $q_1=0,\ q_2=1$, and $q_{i-1}+q_{i+1}=k_i q_i$. Now define the set
\begin{equation*}
K_{e-2}=\{(k_2,\ldots,k_{e-1})\in\mathbb{N}^{e-2}\lsuchthat
	\begin{array}{c}
	\textrm{(i) }[k_2,\ldots,k_{e-1}]\ \textrm{is well defined and yields}\ 0\\
	\textrm{(ii) The corresponding integers}\ q_i\ \textrm{are positive}\end{array}\}.
\end{equation*}
Further, define the set 
\begin{equation*}
K\left(Y_{(n,q)}\right)=\{(k_2,\ldots,k_{e-1})\in K_{e-2} \rsuchthat k_i\leq a_i\}.
\end{equation*}

Each $\underbar{k}\in K\left(Y_{(n,q)}\right)$ determines a fan: $\Sigma_{\underbar{k}}$ is built from the rays generating $\sigma$ and those lying in $\sigma$ which are orthogonal to $w^i/q_i-w^{i-1}/q_{i-1}\in M_\mathbb{R}$ for some $i=3,\ldots,e-1$. Equivalently, the affine lines $\left[\langle\cdot,w^i\rangle=q_i\right]$ form the ``roofs'' of the (possibly degenerate) $\Sigma_{\underbar{k}}$-cones $\tau_i$. The length in the induced lattice of each roof is $(a_i-k_i)q_i$, and this segment lies in height $q_i$.

\begin{thm}
	The P-resolutions of $Y_{(n,q)}$ are in one-to-one correspondence to the elements of $K\left(Y_{(n,q)}\right)$. This correspondence can be realized by the map $\underbar{k}\mapsto \tv\left(\Sigma_{\underbar{k}}\right)$, that is, $\underbar{k}$ corresponds to the toric variety determined by the fan $\Sigma_{\underbar{k}}$.
\begin{proof}
	See \cite{MR1129040} and \cite{MR1652476}.
\end{proof}
\end{thm}

For each $\underbar{k}$, denote by $S_{\underbar{k}}$ the versal base component corresponding to the P-resolution $\tv(\Sigma_{\underbar{k}})$. We will present examples of several fans corresponding to P-resolutions in section \ref{example}..

\begin{remark}
The continued fraction $[1,2,2,\ldots,2,1]=0$ always belongs to $K\left(Y_{(n,q)}\right)$. The P-resolution defined by the corresponding fan is the so-called RDP-resolution of $Y_{(n,q)}$. This corresponds to the Artin component of the versal base space, which always has maximal dimension.
\end{remark}

\section{Milnor Numbers}\label{milnor}
Altmann notes in the introduction of \cite{MR1652476} that T-singularities are exactly those cyclic quotients coming corresponding to a cone $\sigma$ attained by taking the cone over some line segment of integral length $\mu+1$ in lattice height $1$. In such a case, the Milnor number $b_2(F)$ of the $\mathbb{Q}$-Gorenstein smoothing of the singularity is equal to $\mu$. Note that the corresponding P-resolution is simply the identity. On the other hand, Stevens has proven the following general formula:

\begin{prop}\label{milnorformula}
Let $Y$ be a cyclic quotient singularity and let $F$ be the Milnor fiber of a one-parameter smoothing of $Y$ over $S_{\underbar{k}}$. Then $$b_2(F)=\dim T_Y^1-3(e-3)+\#\{2<i<e-1\lsuchthat q_i=1\}+2.$$
\begin{proof}
See lemma 5.3 in \cite{MR1129040}.
\end{proof}
\end{prop}

We will now formulate this in toric terms. Fix a cyclic quotient singularity $Y$ and let $\underbar{k}\in K(Y)$. For any two-dimensional cone $\tau\in\Sigma_{\underbar{k}}^{(2)}$, let $l(\tau)$ and $h(\tau)$ respectively denote the lattice length and height of its roof. We then have the following theorem:

\begin{thm}\label{mymilnor}
	Let $F$ be the Milnor fiber of a one-parameter smoothing of $Y$ over $S_{\underbar{k}}$. Then $$b_2(F)=\bigg(\sum_{\tau\in\Sigma_{\underbar{k}}^{(2)}} l(\tau)/h(\tau)\bigg)-1.$$
\begin{proof}
	The Milnor number $b_2(F)$ can be read from the P-resolution as the sum of the Milnor numbers of the $\mathbb{Q}$-Gorenstein smoothings of each T-singularity plus the total number of exceptional divisors. For any $\tau\in\Sigma_{\underbar{k}}^{(2)}$, the Milnor number of  the $\mathbb{Q}$-Gorenstein smoothing of the T-singularity $\tv(\tau)$ is $l(\tau)/h(\tau)-1$ by Altmann's result. On the other hand, the number of exceptional divisors is simply the number of internal rays in $\Sigma_{\underbar{k}}$, that is, one less than the number of two-dimensional cones $\tau\in\Sigma_{\underbar{k}}^{(2)}$. Combining these two facts yields the above formula.
\end{proof}
\end{thm}
Note that this result generalizes Altmann's result for the special component of a T-singularity.

\section{Component Dimension}\label{formulas} 
We now state the explicit formula for computing the dimension of versal base components for a cyclic quotient singularity $Y$; this comes from Christophersen's and Stevens' description of the versal components.

\begin{prop}\label{dimformula2}
	The dimension of the versal base component $S_{\underbar{k}}$ corresponding to the continued fraction $\underbar{k}\in K(Y)$ can be computed as 
	\begin{equation}\label{eform}
		\dim S_{\underbar{k}} =\#\{2<i<e-1\lsuchthat q_i=1\}+\sum_{i=2}^{e-1}(a_i-k_i) 
	\end{equation}
where the $a_i$, $k_i$, and $q_i$ are as in section \ref{prelim}.
	\begin{proof}
		This formula follows directly from Christophersen's definition of $V_{[\mathbf{k}]}$ in section 2.1.1 of \cite{MR1129026}.\footnote{Note that Christophersen's indices are shifted by one.}
	\end{proof}
\end{prop}

We now translate this formula into toric terms to attain a new dimension formula depending only upon the geometry of the fans of the minimal resolution and the corresponding P-resolution. Let $Y=\tv(\sigma)$ be a cyclic quotient singularity with minimal resolution $\widetilde{Y}=\tv(\Sigma)$ and some P-resolution $\tv(\Sigma_{\underbar{k}})$ corresponding to the versal base component $S_{\underbar{k}}$. Let $v^i$ be as in section \ref{prelim} the generators of the rays in the fan $\Sigma$. For any two-dimensional cone $\tau\in\Sigma_{\underbar{k}}^{(2)}$, let once again $l(\tau)$ and $h(\tau)$ respectively denote the lattice length and height of its roof. Finally, define \begin{equation*}\nu=\sum_{i=1}^r (\det(v^{i-1},v^{i+1})),\end{equation*} that is, $\nu$ is the sum over all $r$ interior rays $v^i$ in $\Sigma$ of the normed volume of the simplex $\conv\{0,v^{i-1},v^{i+1}\}$. This leads to our formula:

	\begin{thm}\label{newform} The dimension of the reduced versal base component $S_{\underbar{k}}$ is given by:
		\begin{equation*} \dim S_{\underbar{k}}=\nu-3r+2\cdot\sum_{\tau\in {\Sigma'}^{(2)}} l(\tau)/h(\tau)-2.\end{equation*}
			\begin{proof}	From corollary 3.18 in \cite{MR608599}, we have $$\dim S_{\underbar{k}}=h^1(\Theta_{\widetilde{Y}})+2b_2(F)-2r,$$
				where $F$ is the Milnor fiber for the component $S_{\underbar{k}}$. Now, $h^1(\Theta_{\widetilde{Y}})=\sum_{i=1}^r (b_i-1)$, see for example \cite{MR0367276}. But each $b_i$ can be computed as $\det(v^{i-1},v^{i+1})$. The desired equation then follows from theorem \ref{mymilnor}.
	\end{proof}	
\end{thm}

When comparing the dimension of two components $S_1$ and $S_2$ with corresponding fans $\Sigma_1$ and $\Sigma_2$ we can even forget about $\nu$ and $r$:
\begin{cor}
The difference in dimension between $S_1$ and $S_2$ is given by:
\begin{equation*}\dim S_1-\dim S_2 = 2\sum_{\tau\in {\Sigma_1}^{(2)}} l(\tau)/h(\tau)-2\sum_{\tau\in {\Sigma_1}^{(2)}} l(\tau)/h(\tau).\end{equation*}
	\begin{proof}
		We express $\dim S_1$ and  $\dim S_2$ using proposition \ref{newform}. The term $\nu-3r-2$ cancels, leaving the desired expression.\end{proof}

\end{cor}

\section{An Example}\label{example}

\begin{figure}
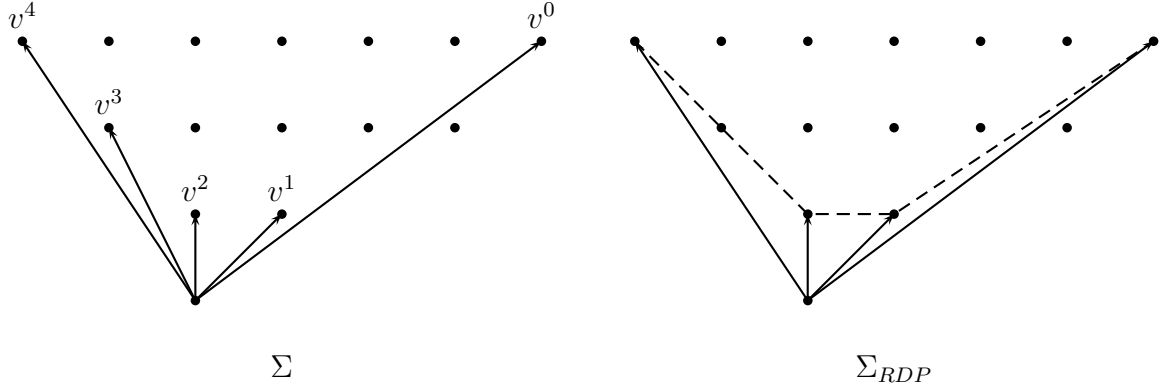

\begin{center}
	$\begin{array}{c@{\qquad\quad} c}
\psset{unit=1.15cm}
\pspicture(-2, 0)(4, 3.5)
\psline{->}(-2,3)
\psline{->}(4,3)
\psdots(0,0)(0,1)(0,2)(-1,2)(1,2)(-2,3)(-1,3)(0,3)(1,3)(1,1)(2,2)(2,3)(3,3)(4,3)(3,2)
\psline{->}(-1,2)
\psline{->}(0,1)
\psline{->}(1,1)
\rput(-2,3.3){$v^4$}
\rput(-1,2.3){$v^3$}
\rput(0,1.3){$v^2$}
\rput(1,1.3){$v^1$}
\rput(4,3.3){$v^0$}

\endpspicture&
\psset{unit=1.15cm}
\pspicture(-2, 0)(4, 3.5)
\psline{->}(-2,3)
\psline{->}(4,3)
\psdots(0,0)(0,1)(0,2)(-1,2)(1,2)(-2,3)(-1,3)(0,3)(1,3)(1,1)(2,2)(2,3)(3,3)(4,3)(3,2)
\psline{->}(0,1)
\psline{->}(1,1)
\psline[linestyle=dashed](-2,3)(0,1)(1,1)(4,3)
\endpspicture\\\\

\Sigma&\Sigma_{RDP}

\end{array}$
\caption{Minimal and RDP resolution fans for $\sigma=\left\langle(-2,3),(4,3)\right\rangle$}\label{fanfig}
\end{center}
\end{figure}
\begin{figure}
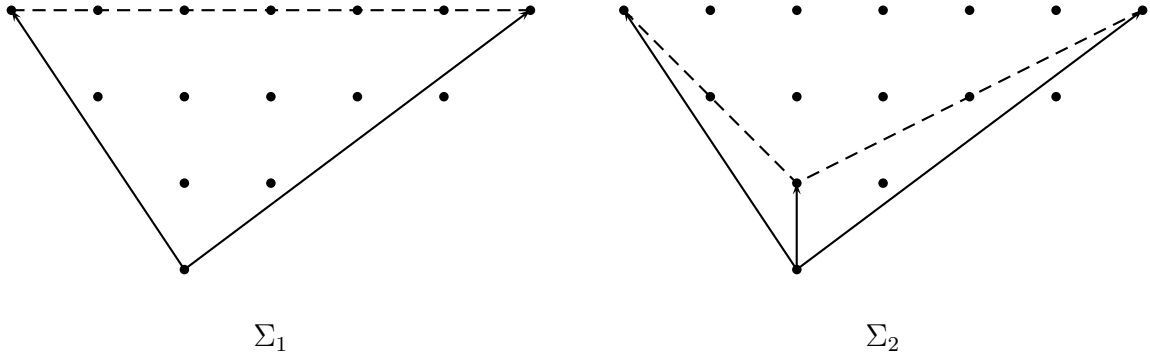

\begin{center}
	$\begin{array}{c@{\qquad\quad} c}
\psset{unit=1.15cm}
\pspicture(-2, 0)(4, 3.5)
\psline{->}(-2,3)
\psline{->}(4,3)
\psdots(0,0)(0,1)(0,2)(-1,2)(1,2)(-2,3)(-1,3)(0,3)(1,3)(1,1)(2,2)(2,3)(3,3)(4,3)(3,2)
\psline[linestyle=dashed](-2,3)(4,3)
\endpspicture&
\psset{unit=1.15cm}
\pspicture(-2, 0)(4, 3.5)
\psline{->}(-2,3)
\psline{->}(4,3)
\psdots(0,0)(0,1)(0,2)(-1,2)(1,2)(-2,3)(-1,3)(0,3)(1,3)(1,1)(2,2)(2,3)(3,3)(4,3)(3,2)
\psline{->}(0,1)
\psline[linestyle=dashed](-2,3)(0,1)(4,3)
\endpspicture\\\\

\Sigma_1&\Sigma_2

\end{array}$
\caption{Further P-resolution fans for $\sigma=\left\langle(-2,3),(4,3)\right\rangle$}\label{fanfig2}
\end{center}
\end{figure}

\begin{ex}
	Consider the cone $\sigma=\langle (-2,3),(4,3)\rangle$ and let $Y=\tv(\sigma)$. It is quite easy to find the minimal resolution fan $\Sigma$; this is done by adding rays through all lattice points on the boundary of $\conv(\sigma\cap N\setminus \{0\})$. Likewise, to get the fan for the RDP resolution, one adds rays through all vertices of $\conv(\sigma\cap N\setminus \{0\})$. These two fans are pictured in figure \ref{fanfig}; note that the dotted lines represent the roofs of two-dimensional cones in $\Sigma_{RDP}$.

	Finding further fans corresponding to P-resolutions is slightly more tricky. Of course, one could use the continued fractions followed by Altmann's construction to get them, but it is also possible without them: every P-resolution is dominated by the resolution corresponding to the fan with rays through all lattice points in $\conv\{0,v^0,v^{r+1}\}$. Thus, by removing these rays and checking for T-singularities and convexity of the roofs, one finds all fans corresponding to P-resolutions. The two additional fans $\Sigma_1$ and $\Sigma_2$ found in this manner are pictured in figure \ref{fanfig2}; note that $Y$ is in fact a $T$-singularity and the fan $\Sigma_1$ corresponds to the space of $\mathbb{Q}$-Gorenstein deformations.

Now, to calculate the desired Milnor numbers and component dimensions, we only need to look at the above pictures. We see that smoothings over the components corresponding to $\Sigma_{RDP}$, $\Sigma_1$, and $\Sigma_2$ have Milnor numbers of $3$, $1$, and $2$, respectively. The number of exceptional divisors $r$ is obviously $3$ and one quickly calculates that $\nu=9$. Thus, the versal base components corresponding to the fans $\Sigma_{RDP}$, $\Sigma_1$, and $\Sigma_2$ have respective dimensions of $6$, $2$, and $4$.

Of course, one could have also calculated that $Y$ corresponds to the chain $(a_2,a_3,a_4,a_5)=(3,3,2,2)$ and that the fans $\Sigma_{RDP}$, $\Sigma_1$, and $\Sigma_2$ correspond, respectively, to the chains $[1,2,2,1]$, $[3,1,2,2]$, and $[2,3,1,2]$. Putting these values into Stevens' formulas yields the same results as with the toric formulas.

\end{ex}
	
\bibliography{p-res_dim}
\address{Mathematisches Institut\\
Freie Universit\"at Berlin\\
Arnimallee 3\\
14195 Berlin, Germany}{nilten@cs.uchicago.edu}
\end{document}